\documentclass[a4paper,11pt]{article}
\RequirePackage{amsthm,amssymb,amsmath}

\def\bB{{\mathbb B}}
\def\bC{{\mathbb C}}
\def\bR{{\mathbb R}}
\def\bP{{\mathbb P}}
\def\bN{{\mathbb N}}

\newtheorem{thm}{Theorem}
\newtheorem{prop}[thm]{Proposition}
\newtheorem{lem}[thm]{Lemma}

\newtheorem{rmk}[thm]{Remark}

\author{Olivier Debarre, Gianluca Pacienza and Mihai P\u aun}
\title{Nondeformability of entire curves in projective hypersurfaces of 
high degree}

\date{\today}
\begin{document}

\maketitle

\begin{quote}
{\small{\bf Abstract.} {In this article, we prove that there does not exist 
a family of entire curves in the universal family of 
hypersurfaces of degree $d\geq 2n$ in the complex projective space 
${\mathbb P}^n$. This can be seen as a weak version of the Kobayashi 
conjecture asserting that a general projective hypersurface of 
high degree is hyperbolic in the sense of Kobayashi.}}
\end{quote}
%

Let $X$ be a hypersurface in the projective space $\bP^n$. The Kobayashi 
conjecture claims that $X$ is {\em hyperbolic}, provided that $X$ 
is general and $d=\deg(X)\geq 2n- 1$. 
By the Brody criterion (\cite{Bro}), the hyperbolicity of $X$ is
equivalent to the fact that every holomorphic map 
$\bC\to X$ is constant.
The conjecture has been proved for $n=3,\ d\geq 21$ and $X$ very general
(\cite{DEG}, \cite{MQ}; see also \cite{Bru1} for an account).
An important and recent progress in the direction of the conjecture 
for all $n\geq 3$ was made by Y.-T. Siu 
in \cite{S}: he obtains a confirmation of the conjecture  
under the assumption   
$d\gg n$.

Consider  the universal family 
${\mathcal X}\subset \bP^n\times \bP^{N_d}$ of hypersurfaces in 
$\bP^n$ with fixed degree $d$ (the number $N_d$ is equal to 
$\binom{n+d}{d} -1$). We will denote by $X_t$ the fiber of 
${\mathcal X}$ over the parameter $t\in \bP^{N_d}$.

\medskip
\noindent 
{\bf Theorem.} {\em Let $U\subset \bP^{N_d}$ be an open set and let
$\Phi:\bC\times U\to {\mathcal X}$ be a holomorphic map
such that $\Phi(\bC\times\{ t\})\subset
X_t$ for all $t\in U$. 
If $d\ge 2n $,   the rank of $\Phi$ cannot be maximal anywhere.}

\medskip
Of course the   theorem
is an immediate consequence of the Kobayashi conjecture (and hence, 
if the degree $d$ is big enough, of the result of Siu). 
So, the theorem above may be rephrased as follows : the Kobayashi conjecture 
may possibly fail only if there is an entire curve on 
a general hypersurface $X$ which is not preserved by a deformation of $X$.

The question above is motivated by the 
``picture'' in the algebraic situation: 
the existence of an algebraic cycle on the 
general member of the family implies its deformation on the nearby fibers. However, 
dealing with transcendental objects (e.g., entire curves) seems to be much more complicated.
For entire curves tangent to a holomorphic foliation of dimension one, 
a substitute for the Hilbert scheme was 
found by M. Brunella in \cite{Bru2}.

Very roughly, the proof goes as follows. First of all, we
consider
the (non-zero) section of the holomorphic bundle
$\Lambda^{1+ N_d}\Phi^*T_{{\mathcal X}}$ given by the Jacobian of $\Phi$
(in fact, for some technical
reasons, 
we will work with a sequence of 
reparame\-trizations of $\Phi$, but we skip this point here, 
to keep the
discussion clear).
In order to use the positivity of the canonical
bundle of the hypersurfaces, we take the wedge product of the previous
section with 
an appropriate family of meromorphic vector fields on ${\mathcal X}$, and
thus  get a
section $\sigma$ of (a twist of) $\Phi^*K^{-1}_{{\mathcal X}}$.
Next, we show that the laplacian of the logarithm of the norm 
of this section dominates a positive multiple of the norm of $\sigma$, and
use negative curvature arguments to derive a contradiction,
as soon as the degree $d$ satisfies the numerical hypothesis of our theorem.

\medskip

\noindent{\bf Proof of the theorem.}
The proof uses two ingredients: the first   is that
the vector bundle $T_{\mathcal X}\otimes p^*{\mathcal O}_{\bP^n}(1)$ 
is generated by its global sections (where $p$ is the   projection 
$\bP^n\times \bP^{N_d}\to \bP^n$). 
The second relies on some negative curvature arguments, very much in the
spirit of the Kobayashi-Ochiai theorem (\cite{KO}).

We recall the following proposition, due to Siu (\cite{S}).

\medskip
\begin{prop}[Siu]\label{siu}
The vector bundle 
$T_{\mathcal X}\otimes p^*{\mathcal O}_{\bP^n}(1)$
is globally generated.
\end{prop}
\medskip

The proof of this proposition is given in \cite{S}; we reproduce it here, 
for the convenience of the reader. 
Observe that the global generation of the restriction 
$T{\mathcal X}_{|X_t}\otimes {\mathcal O}_{X_t}(1)$ of 
the same bundle to a fiber $X_t$ 
has been previously proved by Voisin (\cite{V1}, Prop. 1.1) 
who   deduced from it important results about the algebraic hyperbolicity of a 
(very) general hypersurface (for an account of the subsequent developments of Voisin's approach, see \cite{C}).

\begin{proof}
Consider 
global coordinates $(Z_j)_{0\leq j\leq n}$ 
(resp. $(a_{\alpha})_{\vert \alpha\vert = d})$ on $\bC^{n+ 1}$ (resp. on $\bC^{N_d+ 1}$). The equation of the manifold 
${\mathcal X}$ in $ \bP^n\times \bP^{N_d}$ can be written as
$$ \sum a_\alpha Z^\alpha= 0$$
 where we use here the multi-index notation $\displaystyle 
Z^{\alpha}= \prod Z_j^{\alpha_j}$. Consider the open set 
$U_0= \{Z_0\neq 0\}\times \{a_{d0\dots 0}\neq 0\}$ in $ \bP^n\times \bP^{N_d}$.
For the rest of the proof, we will work on $U_0$, 
with the induced nonhomogeneous
coordinates. 

Consider a multi-index $\alpha\in \bN^d$ and an integer $j$
such that $\alpha_j\geq 1$. On 
the set $U_0$, consider the vector field
$$V_{\alpha, j}=
{\frac{\partial}{\partial a_{\alpha}}}-
z_j{\frac{\partial}{\partial {\hat a}_{ \alpha}}}
$$
 where $z_j= Z_j/Z_0$,     ${\hat a}_k= a_k$  if 
$k\neq j$, and  ${\hat a}_j= \alpha_j- 1$. The vector field 
$V_{\alpha, j}$ is tangent to 
${\mathcal X}_0= {\mathcal X}\cap U_0$, 
as a quick verification shows. On the other hand,
we can extend it to the whole manifold $\mathcal X$ as a meromorphic vector
field and its pole order is equal to $1$. Remark that $V_{\alpha, j}$
is a meromorphic section of the kernel of the differential of the 
first projection  $p_{|{\mathcal X}}: {\mathcal X}\to \bP^n$.

 We also have a ``lifting'' property for the vector fields, as follows.
Consider  a vector field
$$V_0= \sum_{j=1}^nv_j{\frac{\partial}{\partial z_j}}$$
 on $\bC^n$, where 
$\displaystyle v_j= \sum _kv_k^{(j)}z_k+ v_0^{(j)}$ is a polynomial 
of degree at most one in the $z_j$-variables. There exists a vector field
$$V= \sum_{\vert \alpha \vert \leq d}
v_\alpha {\frac{\partial}{\partial a_\alpha}}+ V_0$$
that is tangent to ${\mathcal X}_0$ 
and that extends to the whole manifold
$\mathcal X$ as a holomorphic section of the tangent bundle. 
Indeed, if we want $V$ to be tangent to ${\mathcal X}_0$, the condition to
be satisfied is
$$\sum_\alpha v_\alpha z^\alpha+ \sum_{\alpha, j}a_\alpha v_j
{\frac{\partial z^\alpha}{\partial z_j}}= 0$$
and the complex numbers $v_\alpha$ are simply chosen such that 
the coefficient of the monomial $z^\alpha$ in the above equation
is equal to zero. The extension property is also quickly verified, as well
as the global generation of the bundle 
$T_{\mathcal X}\otimes p^*{\mathcal O}_{\bP^n}(1)$
 by the vector fields already constructed. 
The proposition is thus proved. 
\end{proof}

\vskip 15pt

 Consider  a holomorphic map
$\Phi: \bC\times U\to {\mathcal X}$ 
over the base $U\subset \bP^{N_d}$
as in  the theorem.
We  suppose that $\Phi$ has maximal rank. 
If $d\geq 2n$, we are going to derive a contradiction.

As $U$ is an open set, we can shrink it 
and suppose that it is equal to a polydisc
$\bB(\delta_0)^{N_d}$. 
We will 
consider the following sequence of maps
$$
 \Phi_k: \bB(\delta_0k)^{N_d+ 1}\to {\mathcal X}
$$
 given by 
$\displaystyle \Phi_k(z, \xi_1,\dots,\xi_{N_d})= \Phi (z k^{N_d}, 
\frac{1}{k}\xi_1,\dots,\frac{1}{k}\xi_{N_d})$.
The technical reason for which we need to change the radius of the disc 
will be clear in a moment. 
Notice that the initial map $\Phi_1 =\Phi$ is of maximal rank, 
thus the section
$$J_\Phi(z, \xi)= 
{\frac{\partial \Phi }{\partial z}}\wedge {\frac{\partial \Phi }{\partial \xi_1}}\wedge\dots
\wedge {\frac{\partial \Phi }{\partial \xi_{N_d}}}(z, \xi)\in 
\Lambda ^{1+ N_d}T_{{\mathcal X}, \Phi(z, \xi)}\leqno (1) $$
 of
$\Phi^*\Lambda ^{1+ N_d}T_{\mathcal X}$   is not identically zero. 
Let us assume that $J_\Phi({\underline 0}) $ is nonzero in the corresponding vector space. Remark   that 
$J_{\Phi_k}({\underline 0})= J_\Phi({\underline 0})$, 
for any $k\geq 1$, where 
$J_{\Phi_k}\in \Phi_k^*\Lambda ^{1+ N_d}T_{\mathcal X}$ is the section associated to the map
$\Phi_k$ as indicated in (1). 
  is not identically zero, as a section. 
The positivity of the vector bundle $T_{\mathcal X}$ in the parameter 
space directions
  now comes into the picture: thanks to  Proposition \ref{siu},
we can choose $n-2$ vector fields 
$$
 V_1,\dots,V_{n- 2}\in T_{\mathcal X}\otimes p^*{\mathcal O}_{\bP^n}(1)
$$
such that 
$$
 J_{\Phi_k}({\underline 0})\wedge 
 \Phi_k^*\bigl(V_1\wedge\dots\wedge V_{n-2}\bigr)\neq 0
$$ 
in $\displaystyle K_{\mathcal X}^{-1}\otimes 
p^*{\mathcal O}_{\bP^n}(n-2)_{\Phi_k(0)}$. 

With the vector fields previously chosen, we  consider the following
section
$$
 \sigma_k= J_{\Phi_k}\wedge \Phi_k^*\bigl( V_1\wedge \dots\wedge V_{n-2}\bigr)
$$
of the bundle 
$\Phi_k^*(K_{\mathcal X}^{-1}\otimes p^*{\mathcal O}_{\bP^n}(n-2))\bigr)$ 
over the polydisc.
Its value at the origin is independent of $k$, 
and of course nonzero.
If $q$ is the projection of ${\mathcal X}$ on the parameter space 
$\bP^{N_d}$, then under the 
assumption $d\geq 2n$, the restriction of 
$K_{\mathcal X}\otimes {\mathcal O}_{\bP^n}(2-n))$ to
$\pi_2^{-1}(U)$ is ample (eventually after shrinking once again 
the open subset $U$), hence we can endow this bundle with a metric $h$ of positive curvature.

We now define a sequence of  functions 
$f_k: \bB(\delta_0k)^{N_d+ 1}\to \bR_+$ as follows 
$$\forall w\in \bB(\delta_0k)^{N_d+ 1}\quad
 f_k (w)= 
 \Vert \sigma_k (w) \Vert ^{2/(N_d+1)}_{\Phi_k ^*h^{-1}}.
$$

\begin{rmk}\label{f_k(0)}\upshape
  Notice that, by construction,  there exists a positive number $c$ such that 
for each $k\geq 1$, we have $f_k(0)= c$.
\end{rmk}

We have the following lemma.

\begin{lem}\label{lemma1}
 For each $k\geq 1$, there exists a positive constant $C$ such that 
 we have $\Delta \log f_k\geq Cf_k$ pointwise over 
 the polydisc $\bB(\delta_0k)^{N_d+ 1}$.
\end{lem}
\begin{proof} 
First, remark that by construction, the image of the map $\Phi_k$ 
lies inside $q^{-1}(U)$,
for each $k\geq 1$, so that  
$$
 i\partial{\bar \partial} 
 \log\Vert \sigma_k \Vert ^2_{\Phi_k^*h^{-1}}\geq 
 \Phi_k^*\Theta_h\bigl(K_{\mathcal X}\otimes 
 p^*{\mathcal O}_{\bP^n}(2-n)\bigr). 
 \leqno (2)
$$
In the inequality (2), take the trace 
with respect to the flat
metric on the polydisc. We get
\begin{eqnarray*}
 \Delta \log\Vert \sigma_k \Vert ^2_{\Phi_k^*h^{-1}}
 & \geq &
 C\Bigl(\Big \Vert {\frac{\partial \Phi_k}{\partial z}}\Big \Vert ^2_\omega
 + \sum_{j=1}^{N_d}
\Big \Vert {\frac{\partial \Phi_k} {\partial \xi_j}}\Big \Vert ^2_\omega\Bigr)\\
& \geq & 
 C\Vert J_{\Phi_k}\Vert ^{2/(1+N_d)}_{\Lambda^{1+ N_d}\omega}\\
& \geq & 
 C\Vert \sigma_k \Vert _{\Phi_k^*h^{-1}}^{2/(1+N_d)} 
\end{eqnarray*}
The constant $C$ in the previous sequence of inequalities 
varies from one line to another, but we still denote it by $C$ 
as it is independent of $k$.
The above relations are obtained using the vector inequalities
$$
 \Vert W_1 \wedge \dots \wedge W_s \Vert
 \leq \Vert W_1 \Vert \dots \Vert W_s \Vert
 \leq s^{-s} (\Vert W_1 \Vert + 
 \dots + \Vert W_s \Vert)^s.
$$
So the lemma is proved.
\end{proof}

Using the previous lemma we will prove a result whose proof is very close 
to that of the Ahlfors-Schwarz lemma. 

\begin{prop}\label{tozero}
 For each $k\geq 1$ we have $f_k(0)\leq Ck^{-2}$. 
 In particular, as $k\to \infty$, we have $f_k(0)\to 0$. 
\end{prop} 
\begin{proof}
Consider the 
volume form of the Poincar\'e metric on the polydisc
$$\psi_k= \frac{1}{\Bigl(1- \frac{\vert z\vert ^2}{\delta_0^2k^2}\Bigr)^2}
\prod_{j=1}^{N_d}\frac{1}{\Bigl(1- \frac{\vert \xi_j\vert ^2}{\delta_0^2k^2}\Bigr)^2}$$
A quick computation shows that 
$$
 i\partial{\bar \partial} \log \psi_k= (\delta_0 k)^{-2}
 \Biggl(\frac{idz\wedge d\overline z}{\Bigl(1- \frac{\vert z\vert ^2}{\delta_0^2k^2}\Bigr)^2}+
 \sum_{j= 1}^{N_d} \frac{id\xi_j\wedge d\overline{\xi _j}}{\Bigl(1-\frac{\vert \xi_j\vert ^2} {\delta_0^2k^2}\Bigr)^2}\Biggr)
$$
so if we take the trace of this equality with respect to the flat metric, we  get
$$\Delta \log \psi_k\leq Ck^{-2}\psi_k.\leqno (3)$$
Remark that the previous inequality can be obtained precisely because we have the
same radius $\delta_0k$ for the components of the 
polydisc which is the domain of
$\psi_k$. This is why we had to reparametrize 
our map $\Phi$ from the very beginning.

Consider the function $\displaystyle (z, \xi)\mapsto \frac{f_k(z, \xi)}{\psi_k(z, \xi)}$. Its maximum cannot be achieved at a boundary point of the domain, since 
$\psi_k$ goes to infinity as $(z, \xi)$ goes to the boundary.
So at the maximum point 
$(z_0, \xi_0)$, we have 
$$\Delta \log f_k/\psi_k\leq 0.\leqno (4)$$
 This inequality, combined with Lemma \ref{lemma1} 
and (3), gives
$$f_k(z_0, \xi_0)\leq Ck^{-2}\psi_k(z_0, \xi_0).\leqno (5)$$
Since the relation (5) is verified at the maximum point of the quotient,
it follows that the same is true at an arbitrary point, so we get 
$$f_k(z, \xi)\leq Ck^{-2}\psi_k(z, \xi).\leqno (6)$$
To finish the proof, it is sufficient to 
write the inequality (6) at the origin.\end{proof}

Since  Proposition \ref{tozero} 
and  Remark \ref{f_k(0)} contradict each other, the theorem is proved.
\hfill $\Box$

\noindent
{\tt debarre@math.u-strasbg.fr},
{\tt pacienza@math.u-strasbg.fr}\\
Institut de Recherche Math\'ematique Avanc\'ee\\
Universit\'e L. Pasteur et CNRS
\\7, rue Ren\'e Descartes, 67084 Strasbourg C\'edex - France

\medskip

\noindent
{\tt Mihai.Paun@iecn.u-nancy.fr}\\
Institut \'Elie Cartan\\ 
Universit\'e Henri Poincar\'e Nancy 1\\
B.P. 239, F-54506 Vand\oe uvre-l\`es-Nancy C\'edex - France.

\end{document}